\documentclass[times,doublespace]{nmeauth}

\usepackage{moreverb}

\usepackage{bm}

\usepackage[colorlinks,bookmarksopen,bookmarksnumbered,citecolor=red,urlcolor=red]{hyperref}

\newcommand\BibTeX{{\rmfamily B\kern-.05em \textsc{i\kern-.025em b}\kern-.08em
T\kern-.1667em\lower.7ex\hbox{E}\kern-.125emX}}

\def\volumeyear{2010}

\begin{document}

\runningheads{M.~U.~Altaf et. al.}{Downscaling the 2D B\'{e}nard Convection Equations Using CDA}

\title{Downscaling the 2D B\'{e}nard Convection Equations Using Continuous Data Assimilation\footnotemark[2]}

\author{M. U. ~Altaf\affil{1}, E. S. ~Titi\affil{2}, O. M. ~Knio\affil{1}, L. ~Zhao\affil{3}, M. F. ~McCabe\affil{1} and I. ~Hoteit\affil{1}\corrauth}

\address{\affilnum{1}King Abdullah University of Science and Technology, Saudi Arabia\break
\affilnum{2}Department of Mathematics, Texas A$\&$M University, College Station, TX 77843, USA\break
\address<\affilnum{3}Georgia Institute of Technology, Atlanta, GA, USA}

\corraddr{King Abdullah University of Science and Technology, Saudi Arabia. E-mail: ibrahim.hoteit@kaust.edu.sa}

\begin{abstract}
We consider a recently introduced continuous data assimilation (CDA) approach for downscaling a coarse resolution configuration of the 2D B\'{e}nard convection equations into a finer grid. In this CDA, a nudging term, estimated as the misfit between some interpolants of the assimilated coarse grid measurements and the fine grid model solution, is added to the model equations to constrain the model. The main contribution of this study is a performance analysis of CDA for downscaling measurements of temperature and velocity. These measurements are assimilated either separately or simultaneously and the results are compared against those resulting from the standard point-to-point nudging approach (NA). Our numerical results suggest that the CDA solution outperforms that of NA, always converging to the true solution when the velocity is assimilated as has been theoretically proven. Assimilation of temperature measurements only may not always recover the true state as demonstrated in the case study. Various runs are conducted to evaluate the sensitivity of CDA to noise in the measurements, the size and the time frequency of the measured grid, suggesting a more robust behaviour of CDA compared to NA.

\end{abstract}


\maketitle



\section{Introduction}
\vspace{-2pt}
An important use of present day numerical models is to provide accurate predictions of physical phenomena from a known present state and some prior information. In some cases, the output from these numerical models is often too spatially coarse to be used directly to assess particular local phenomena (e.g. climate models). In order to improve the representation of local scale processes, downscaling techniques can be applied. Downscaling methods are needed to obtain local scale weather and climate, based on regional scale atmospheric variables that are provided by global circulation models (GCMs). Two basic approaches are followed in downscaling: $i)$ statistical downscaling and $ii)$ dynamical downscaling.

In statistical downscaling (SD), an empirical relationship is established from observations between large scale variables, and a local variable at a particular station. The relationship is then used to estimate the local variables from the coarse global fields \cite{Timbal2003, Hewitson2006, Gutzler2011, Jha2015}. A review of SD methods can be found in \cite{Guitrez2013}. Dynamical downscaling (DD) on the other hand, fits regional climate models (RCMs) into outputs from GCMs \cite{Mcgregor1997}. The nudging algorithm is one of the standard methods in DD, which basically forces RCM simulations toward large-scale driving data point-to-point. A number of studies have examined the performance of  nudging in reanalysis-driven RCM simulations (e.g. \cite{Liu2012, Feser2012}). These RCM simulations show better temporal variations of precipitation relative to those without nudging \cite{Lo2008}. Various studies demonstrated that the overall performance of SD and DD were comparable in reproducing the present day climate for the respective regions \cite{Wilby1997, Murphy1999}. It is generally more straightforward and physically consistent to use the dynamical information, rather than a statistical model, when available to downscale global fields.  

The study of different geophysical processes is based on both modelling and observations of the dynamical model state. None of the two disciplines has proved to be more useful. Instead both disciplines take advantage of each other: Observational data usually serve to validate numerical models, and observations are interpolated in space and time using the models. Data assimilation (DA) is the technique that merges these two sources into a single product. The basic idea of DA is to regularly update a numerical model with observational data \cite{Bennet1992}. Given a sufficiently realistic numerical model and reliable observations, the joint product provides an estimate of the state of the ocean or atmosphere that should be better than the respective single sources. Apart from the quality of the model or the observations, the success of DA depends on how well the two sources are combined. 


Continuous data assimilation (CDA) methods incorporate observational data directly into the model equations as the model is being integrated in time \cite{Charney1969, Daley1991}. One way to achieve this is by introducing low-Fourier-mode observables into the model equations for the evolution of high Fourier modes \cite{Henshaw2003, Olson2003, Olson2009, Korn2009, HOTi}. Recently \cite{Azouani2014} introduced a new approach to CDA based on ideas from control theory. Instead of inserting the measurements directly into the model, in this CDA approach, a nudging term, estimated as the misfit between some interpolants of the assimilated coarse grid measurements and the fine grid model solution, is added to the model equations to constrain the model large scale variability by available measurements. This new CDA approach was designed for general dissipative dynamical systems and has been successfully analyzed in different scenarios \cite{Bessaih2015, Farhat2015a, Farhat2015b}, and tested numerically in \cite{GOTi}. 

In this paper, the new CDA  method is examined for downscaling a coarse scale solution of the non-linear 2D B\'{e}nard convection problem (a Boussinesq system between two solid plates, heated from the bottom and cooled from the top). The B\'{e}nard convection equations solve the reduced form of Navier-Stokes equations and are one of the classical problems in the heat transfer literature \cite{Aswatha2012}. The coarse mesh data is generated for temperature and velocity. The CDA algorithm is then applied to downscale these variables, assimilating temperature observations, velocity observations or both. The CDA algorithm in \cite{Azouani2014} has been theoretically investigated in \cite{Farhat2015} for the 2D B\'{e}nard convection through velocity measurements alone,  with a finding that velocity measurements alone are enough to construct the approximate solutions for both velocity and temperature. Earlier study by Charney et. al. \cite{Charney1969} claimed that the temperature measurements are enough to determine all other variables in the atmosphere. Ghil et. al in \cite{Ghil1977, Ghil1978} showed that this was only true for some simple models, but gave other numerical tests showing that this is not always true. Here, we perform a number of numerical experiments with CDA assimilating velocity and temperature and compare the results with the standard nudging approach \cite{Hoke1976}. We also study the performance of CDA with respect to perturbations in the assimilated observations, time frequency of observations, the number of observation points, and different choices of interpolants. Our results suggest that the CDA method is a fast and accurate approach with favorable robustness properties for dynamical downscaling of the 2D B\'{e}nard convection equations. 

In the next section a description of the 2D B\'{e}nard convection numerical model is given, after which we briefly describe the CDA technique in Sections \ref{sec.cda}. The numerical results obtained with the 2D B\'{e}nard convection model are presented in Section \ref{sec.exp} while  conclusions are provided in Section \ref{sec.con}.

\section{B\'{E}nard convection problem}
Known as the B\'{e}nard convection problem, the thermally driven rectangular cavity with adiabatic top and bottom walls is one of the classical problems in the heat transfer literature \cite{Aswatha2012}. It is also one of the most popular test-problems for comparing numerical algorithms. The model equations are casted in dimensionless form using an appropriate scaling 
of dynamical variables.  Using tildes to denote dimensional quantities, we use the cavity height, $\tilde{H}$ as reference lengthscale, $V_r \equiv (\tilde{\kappa} / \tilde{H})Ra^{0.5}$ as
reference velocity, and the background temperature difference $\tilde{T}_2 - \tilde{T}_1$ as characteristic temperature.  $Ra$ denotes the Rayleigh number, whereas $\tilde{\kappa}$ 
is the thermal conductivity.  With the current choice of scaling parameters, the normalized momentum, energy and mass conservation equations are:
\begin{equation} 
\frac{\partial \bm{u} }{\partial t} + (\bm{u} \cdot \nabla)\bm{u} +\nabla p =  \frac{Pr}{\sqrt{Ra}} \nabla^2 u + Pr\Theta \bm{e}_{2},
\end{equation}
\begin{equation} 
\frac{\partial \Theta}{\partial t} + (\bm{u} \cdot \nabla)\Theta  = \frac{1}{\sqrt{Ra}} \nabla^2 \Theta + \bm{u} \cdot \bm{e}_{2},
\end{equation}
\begin{equation} 
\nabla \cdot u = 0 .
\end{equation}
where $\bm{u}=(u,v)$ denotes the velocity vector, $p$ is pressure, $t$ is time,  $\Theta \equiv (T- 0.5(T_1+T_2))/(T_2-T_1)$ is the temperature anomaly, and $Pr$ is the Prandtl number.  
The system is completed by specifying the initial conditions:
\begin{equation} 
\bm{u}(0;(x,y)) = \bm{u}_{0}(x,y), 
\end{equation}
and 
\begin{equation} 
\Theta(0;(x,y)) = \Theta_{0}(x,y),
\end{equation}
the boundary conditions and the top and bottom boundaries,
\begin{center}
$\bm{u} = 0 \, \, $at $\, \,  y = 0 \,$ and $\,  y = 1$,
\end{center}
and
\begin{center}
$\Theta = 0 \, \,$at $\, \,  y = 0 \,$ and $\,  y = 1$,
\end{center}
and periodicity conditions in $x$ for velocity, temperature, and pressure.

\subsection{Numerical Solution}

The governing equations are simulated using a finite difference methodology.  We use 
a staggered, marker-and-cell (MAC) grid to discretize field variables.  Specifically, a 
Cartesian grid is used, with $n_x$ and $n_y$ cells in the $x$ and $y$ directions respectively.
The corresponding cell sizes are $\Delta x = 1/n_x$ and $\Delta y = 1/n_y$.  Temperature
and pressure are discretized on cell centers, whereas the velocity components $u$ and
$v$ are discretized on vertical and horizontal cell faces respectively.  Gradients are
approximated using conservative, second-order, centered differences.

The system of discretized equations is integrated in time using a pressure
projection scheme.  The pressure-free momentum and energy equations are first integrated
in time using a mixed integration approach.  The second-order Adams-Bashforth scheme is
used for the convective terms, diffusion terms are treated using the Euler backward scheme,
whereas the remaining (linear) terms are handled using a first-order explicit Euler scheme.
A pressure projection step is then applied in order to satisfy the continuity equation.

\section{Downscaling Approaches}\label{sec.cda}

\subsection{Grid Nudging Method}


Suppose the time evolution of $(\bm{u}, \Theta)$ is governed by a given system of the form (1) -- (5), mentioned above, where the initial data is not known. We construct an increasingly accurate initial condition based on coarse grid data from which predictions of $(\bm{u}(t), \Theta(t))$ can be made. In the standard NA, the  interior solution of the model, denoted as $(\bm{u}(t), \Theta(t))$, is brought closer to, or ``nudged'' towards the observations, $\left\{ \bm{u}^o_j, \Theta^o_j \right\}$, $j=1,\ldots,N$, where $N$ is the number of observation points located at $\bm{x}_j$, namely by adding an adjustment or restoring term.  The velocity nudging term is expressed as: $\sum_{j=1}^N \alpha_u (\bm{u}^o_j - \bm{u}(\bm{x_j})) \delta(\bm{x}-\bm{x}_j)$, where the ``nudging coefficient''  $\alpha_u$ is tuned to minimize discrepancy; a similar nudging is used for temperature.  The updated equations will then be of the form:
\begin{equation} 
\frac{\partial \bm{u}}{\partial t} - \frac{Pr}{\sqrt{Ra}} \nabla^2 \bm{u} + (\bm{u} \cdot \nabla)\bm{u} +\nabla p = \Theta \bm{e}_{2} +  \sum_{j=1}^N \alpha_u (\bm{u}^o_j - \bm{u}(\bm{x_j})) \delta(\bm{x}-\bm{x}_j) ,
\end{equation}
\begin{equation} 
\frac{\partial \Theta}{\partial t} -\frac{1}{\sqrt{Ra}} \nabla^2 \Theta + (\bm{u} \cdot \nabla)\Theta - \bm{u} \cdot \bm{e}_{2} = \sum_{j=1}^N \alpha_\Theta (\Theta^o_j - \Theta(\bm{x_j})) \delta(\bm{x}-\bm{x}_j) ,
\end{equation}
\begin{equation} 
\nabla \cdot \bm{u} = 0,
\end{equation}
with initial conditions:
\begin{equation} 
\bm{u}(0;(x,y)) = 0,
\end{equation}
\begin{equation} 
\Theta(0;(x,y)) = 0.
\end{equation}
We call the solution (1) -- (5) as reference solution, and the solution of (21) -- (25) as estimated solution.

\subsection{Continuous Data Assimilation}

In the continuous data assimilation (CDA) approach, interpolants are applied on both model outputs and coarse data before nudging, allowing to only constrain the large scale flow of the model as shown by \cite{Azouani2014}. Let $I_h(\phi(x))$ be an interpolation operator of a function $\phi(x)$. For instance one can take

\begin{equation} 
I_h(\phi)(x) = \sum_{k=1}^{N_h}\phi(x_k)\chi_{Q_k}(x),
\end{equation}
where $Q_j$ are disjoint subsets such that diam$(Q_j)\leq h, \, \, \bigcup_{j=1}^{N_h}Q_{j} = \Omega , \, \, x_{j}\in Q_{j}, \, \, \chi_E$ is the characteristic function of the set E, and $\phi$ is a suitable interpolant.

Consider the new system of equations:
\begin{equation} 
\frac{\partial \bm{u}}{\partial t} - \frac{Pr}{\sqrt{Ra}} \nabla^2 \bm{u} + (\bm{u} \cdot \nabla)\bm{u} +\nabla p = \Theta \bm{e}_{2} - \mu_{u}(I_h(\bm{u}^o) - I_h(\bm{u})),
\end{equation}
\begin{equation} 
\frac{\partial \Theta}{\partial t} -\frac{1}{\sqrt{Ra}} \nabla^2 \Theta + (\bm{u} \cdot \nabla)\Theta - \bm{u} \cdot \bm{e}_{2} = -\mu_{\theta}(I_h(\Theta^o) - I_h(\Theta)),
\end{equation}
\begin{equation} 
\nabla \cdot \bm{u} = 0,
\end{equation}
with initial conditions:
\begin{equation} 
\bm{u}(0;(x,y)) = 0,
\end{equation}
\begin{equation} 
\Theta(0;(x,y)) = 0.
\end{equation}

Let's denote the estimated solution of (27) -- (31) as $V(t) = (\bm{u}(t), \Theta(t))$. To accurately predict $U(t)$ on the interval $[t_1, \, \, t_1 + T]$, it is sufficient to have coarse observational data $I_h(U^o(t))$ accumulated over an interval of time $[0, \, \, t_1]$ linearly proportional to $T$ in the immediate past \cite{Azouani2014}.

In particular, suppose it is desired to predict $U(t)$ with accuracy $\epsilon > 0$ on the interval $[t_1, \, \, t_1 + T]$, where $t_1$ is the present time and $T > 0$ determines how far into the future one wants to predict. Let $h$ be small enough and $\mu_u$ be large enough, then there exist constants $C$ and $\beta$ such that
\begin{equation} 
\|U(t) - V(t)\|_{H_1(\Omega)} \leq Ce^{-\beta t}, \, \, t\geq 0.
\end{equation}

One then uses $V(t_1)$ as the initial condition to make future predictions. Let $W$ be a solution of (1) -- (5) with initial condition $W(t_1) = V(t_1)$. Known results on continuous dependence on initial conditions imply the existence of $\gamma > 0$ such that
\begin{equation} 
\|W(t) - U(t)\|_{L_2(\Omega)} \leq \|W(t_1) - U(t_1)\|_{L_2(\Omega)} e^{\gamma (t-t_1)}, \, \, t\geq t_1.
\end{equation}
Therefore,
\begin{equation} 
\|W(t) - U(t)\|_{L_2(\Omega)} \leq Ce^{-\beta t_1 + \gamma T}, \, \, t \in [t_1, \, \, t_1 + T],
\end{equation}
provided $\beta t_1 \geq \gamma T + \ln(C/\epsilon)$. Thus $W(t)$ predicts $U(t)$ with accuracy $\epsilon$ on the interval           $[t_1, t_1 + T]$.


\section{Numerical Experiments}\label{sec.exp}

We implemented the 2D B\'{e}nard Convection Model in Matlab, as described in section $(2)$. We consider the box domain $[0, 2] \times [0, 1]$ uniformly discretized into $200 \times 100$ computational cells. The time step in our simulations is equal to $0.01$ and a total of $3000$ integration steps is performed. Our goal is to reconstruct, as accurately as possible the solution ($\bm{u}(t)$; $\Theta(t)$) of the B\'{e}nard convection problem on a high resolution grid, from observations available at a much coarser grid using CDA. The results obtained using CDA are then compared with predictions of standard NA. We study the sensitivity of CDA to different choices of the interpolation operator $I_h$, and investigate the CDA robustness to observational noise and spatial distribution of observations. We also provide an example in which temperature data alone is not enough to estimate the state of the B\'{e}nard equations (i.e., both temperature and velocities).

\begin{figure}[htbp]
\centering
	\includegraphics[height=150mm,width=150mm]{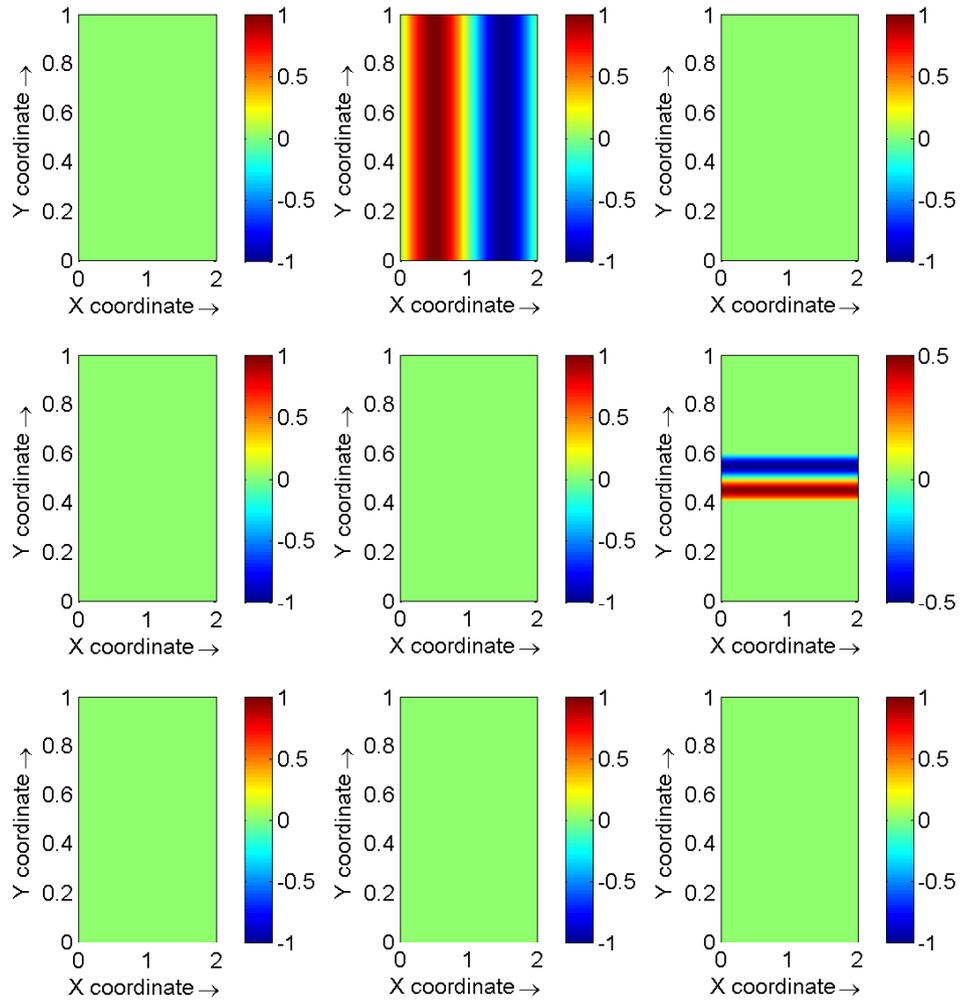}
\vskip 100pt 
	\caption{Initial conditions for the background solution (left panel), reference solution (middle panel), and initial shear flow (right panel). Top: Temperature, Middle: $u$-velocity and Bottom: $v$-velocity}
	\label{fig1}
\end{figure}

A reference model run was conducted with the initial conditions for temperature set as $sin(2 \pi x+2000 \pi y)$, and the velocity initialized from rest (see Figure \ref{fig1}). The outputs of this reference run are used to obtain the observations. The comparison of the two schemes (CDA and NA) is evaluated in terms of relative root-means-square errors (RRMSE), computed as the ratio between the $L_2$-norm of the difference between the reference and model solutions and the $L_2$-norm of the reference solution after each assimilation step.

\subsection{Sensitivity with respect to data grid resolution}\label{sengrid}

Table \ref{tab2} outlines the three values of relaxation parameters ($\mu_\theta$ and $\mu_u$) and nudging operators ($\alpha_\theta$ and $\alpha_u$) tested in the assimilation experiments with CDA and NA respectively referred to as large, medium and small.  We conducted three experiments assimilating $(i)$ velocity, $(ii)$ temperature, and $(iii)$ both velocity and temperature. Three observational scenarios are studied: $1$ out of $20$ grid points of the high resolution grid is assumed observed, for a total of $10 \times 5$ observed points, $1$ out of $10$ grid points, for a total of $20 \times 10$ observed points, and $1$ out of $20$ grid points, for a total of $40 \times 20$ observed points. Data were assimilated at every time step. We started the simulation from zero initial conditions for both temperature and velocity.

\begin{figure}[htbp]
\centering
	\includegraphics[height=150mm,width=150mm]{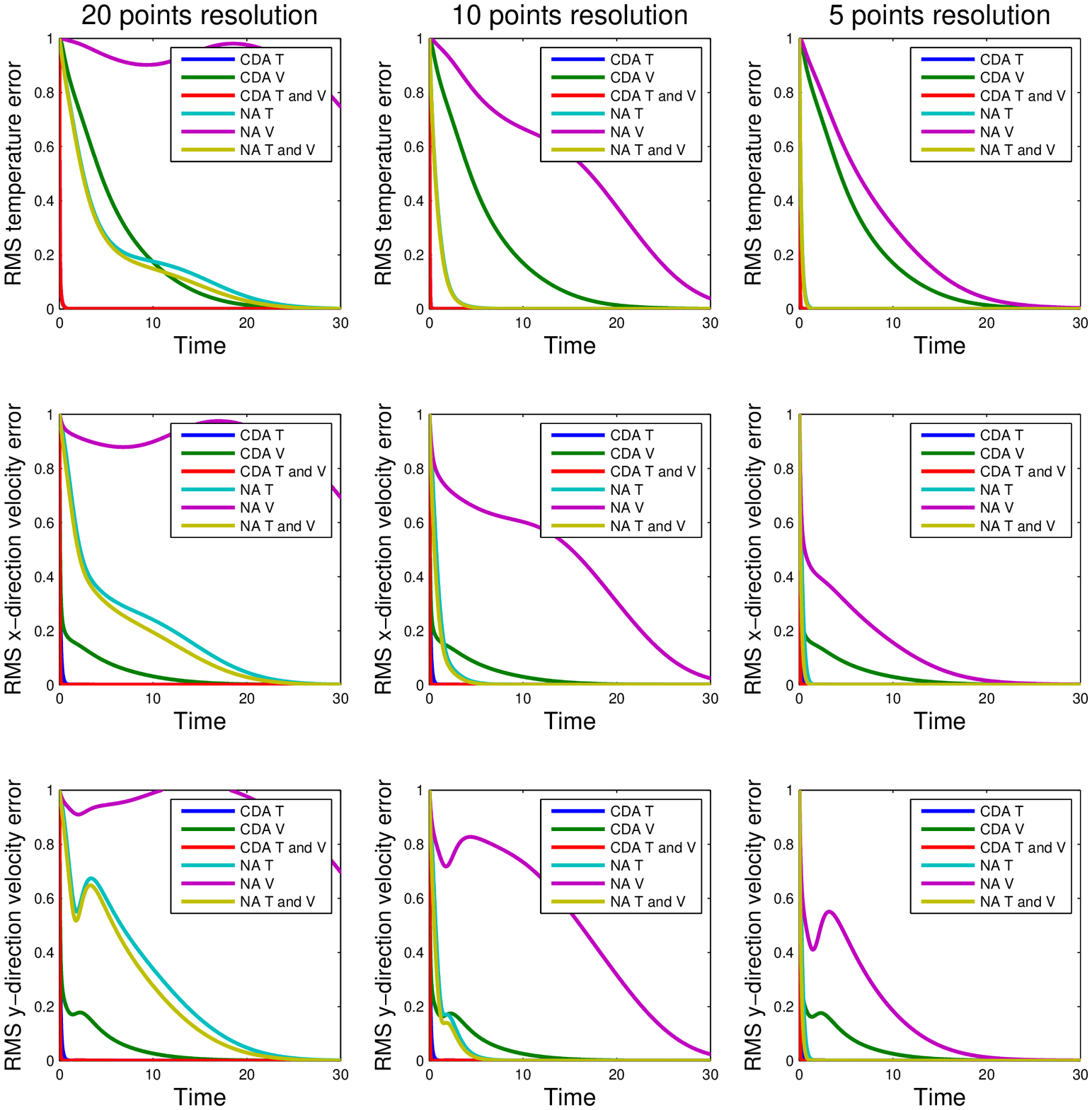}
\vskip 100pt 
	\caption{CDA  vs NA RRMSE results (CDA - $\mu_\theta=0.10, \mu_V = 0.10$, NA - $\alpha_\theta
	=1.50, \alpha_V = 1.50$). The left panel shows 20 point resolution, the middle panel shows 10 point resolution, and  5 point resolution on the right panel. The blue curve (CDA T) is overlapped by the red line (CDA T and V).}
	\label{fig2}
\end{figure}

\begin{figure}[htbp]
\centering
	\includegraphics[height=150mm,width=150mm]{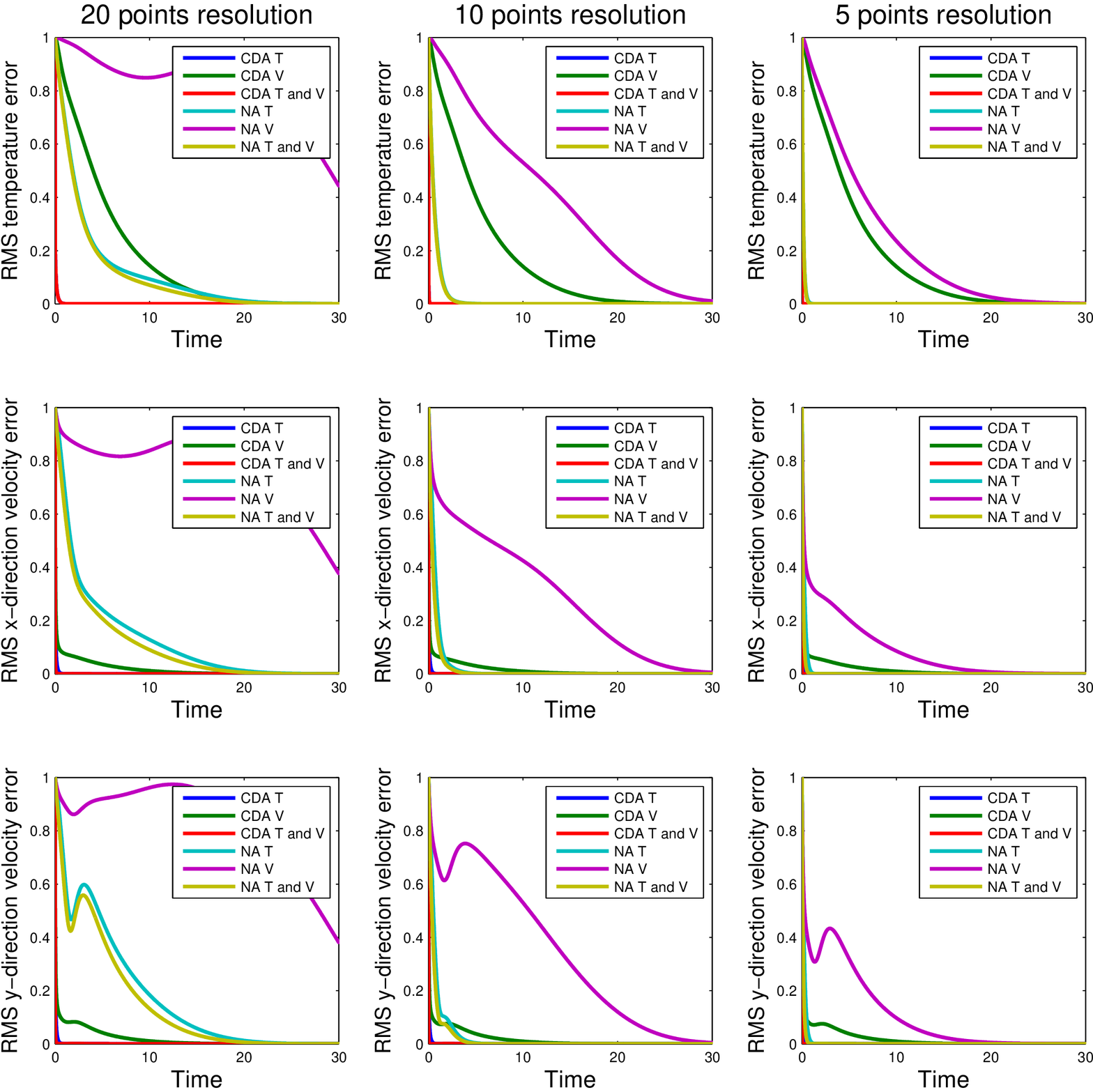}
\vskip 100pt 
	\caption{CDA  vs NA RMSE results (CDA - $\mu_\theta=0.50, \mu_V = 0.50$, NA - $\alpha_\theta=2.50, \alpha_V = 2.50$). The left panel shows 20 point resolution, the middle panel shows 10 point resolution, and  5 point resolution on the right panel. The blue curve (CDA T) is overlapped by the red line (CDA T and V).}
\label{fig3}
\end{figure}

\begin{figure}[htbp]
\centering
	\includegraphics[height=150mm,width=150mm]{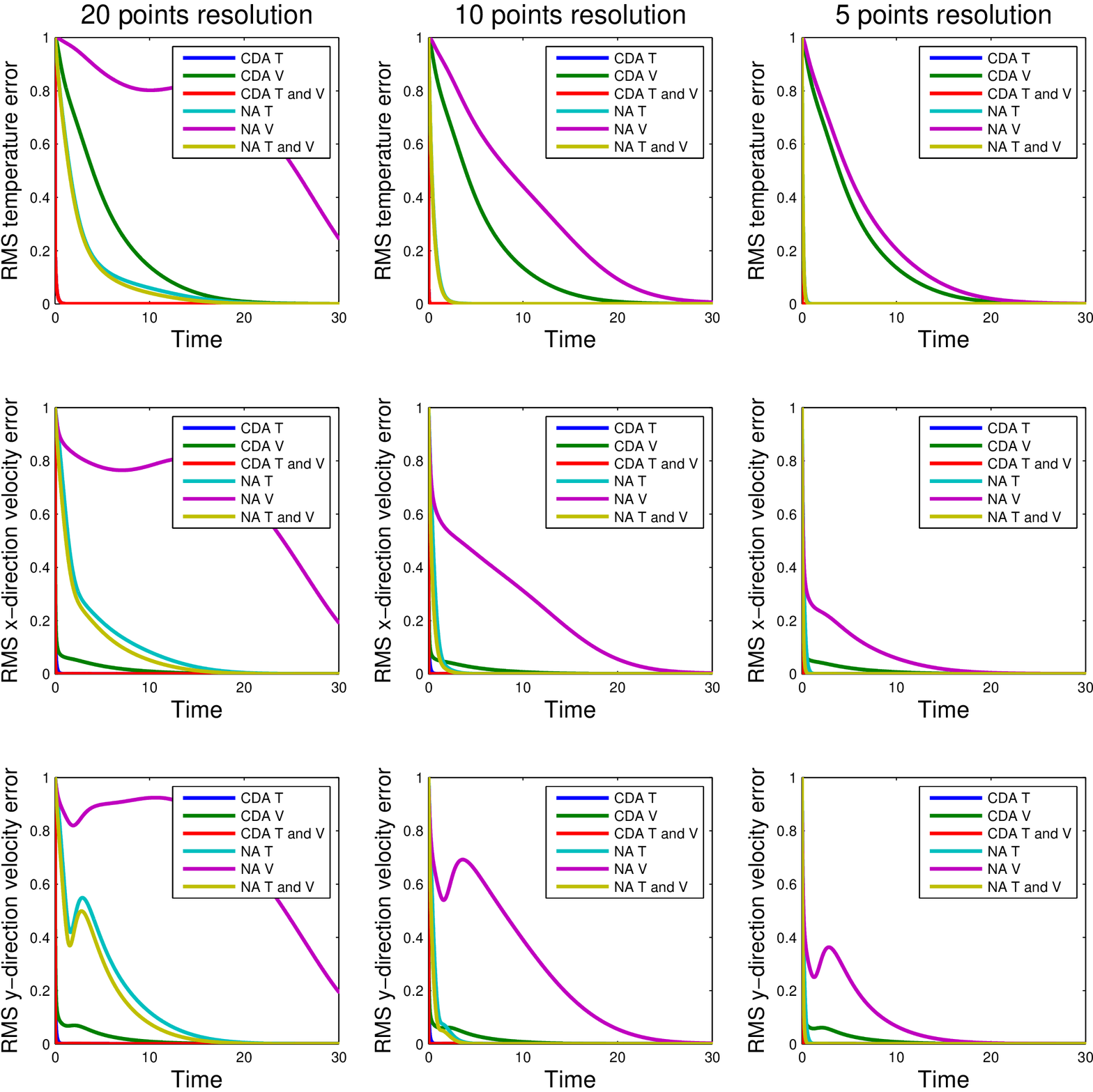}
\vskip 100pt 
	\caption{CDA  vs NA RMSE results (CDA - $\mu_\theta=1.0, \mu_V = 1.0$, NA - $\alpha_\theta=3.50, \alpha_V = 3.50$). The left panel shows 20 point resolution, the middle panel shows 10 point resolution, and  5 point resolution on the right panel. The blue curve (CDA T) is overlapped by the red line (CDA T and V). }
\label{fig4}
\end{figure}

Figures \ref{fig2} - \ref{fig4} plot the time evolution of the  RRMSE with respect to the norm of the reference solution as they result from CDA and NA assimilating temperature, velocities and both temperature and velocities. The results from three different scenarios and using different nudging coefficients suggest that CDA performs very well, showing significant improvements over NA estimates. The convergence of CDA is clearly faster than NA as depicted in Figure \ref{fig2}. When assimilating either temperature (i.e., $\mu_u = 0$) or both velocity and temperature, an  exponential convergence is obtained for both velocity and temperature for all three choices of grid resolution. This is not true for NA, which achieves an exponential convergence only for the case when dense observation network is used (i.e., 5 point grid resolution). 

\begin{figure}[htbp]
\centering
	\includegraphics[height=150mm,width=150mm]{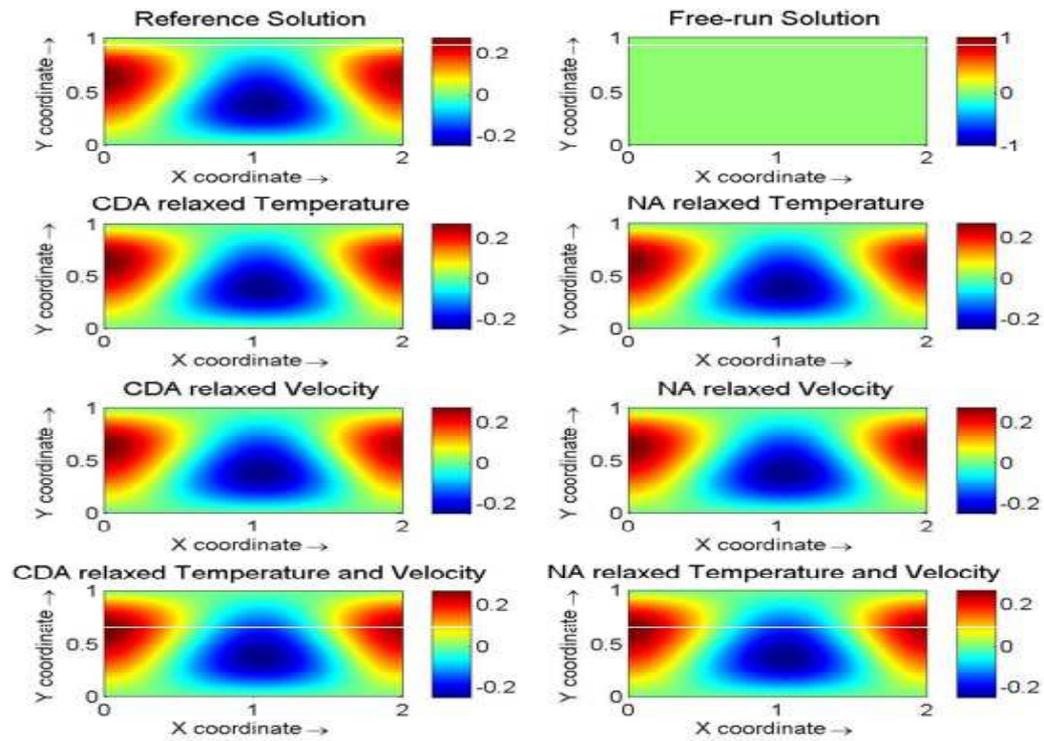}
\vskip 20pt 
	\caption{Temperature field at time ($t = 30$) and assimilation of data every 10 point grid resolution as it results from CDA and NA with assimilation of different variables (CDA - $\mu_\theta=0.10, \mu_V = 0.10$, NA - $\alpha_\theta=1.50, \alpha_V = 1.50$). The simulation starts with zero initial conditions (\ref{fig2} - left panel).}
\label{fig5}
\end{figure}

\begin{figure}[htbp]
\centering
	\includegraphics[height=150mm,width=150mm]{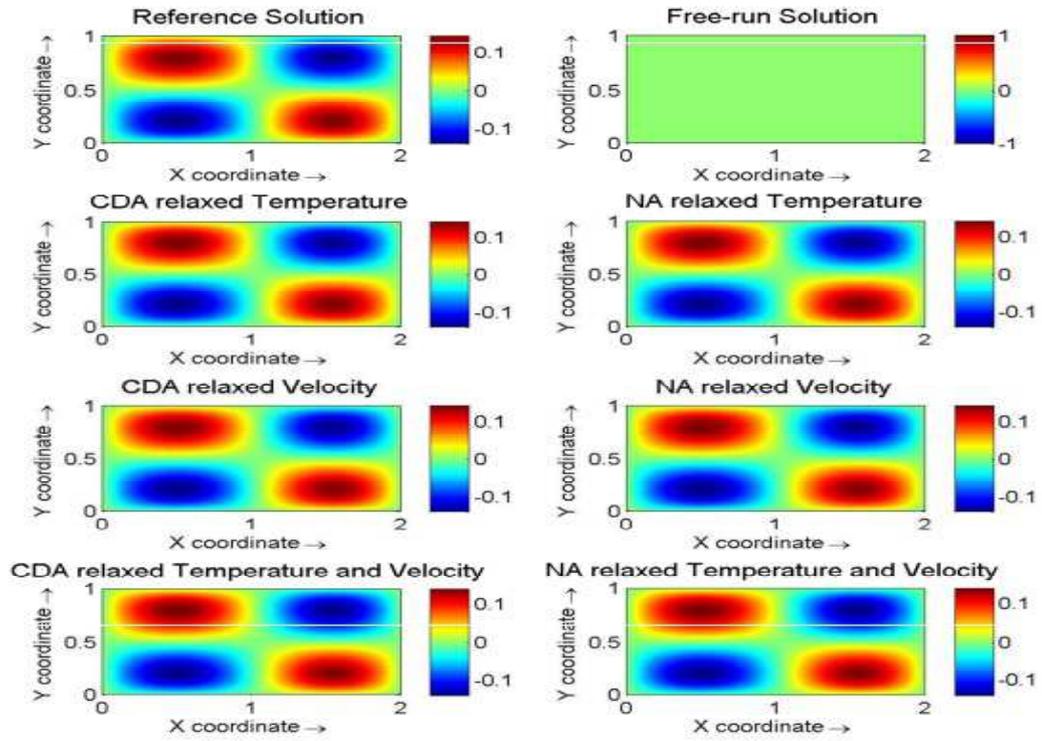}
\vskip 20pt 
	\caption{$x$-direction of velocity field at time ($t = 30$) and assimilation of data every 10 point grid resolution as it results from CDA and NA with assimilation of different variables (CDA - $\mu_\theta=0.10, \mu_V = 0.10$, NA - $\alpha_\theta=1.50, \alpha_V = 1.50$). The simulation starts from rest, i.e. $u=0$ (\ref{fig2} - left panel).}
\label{fig6}
\end{figure}

\begin{figure}[htbp]
\centering
	\includegraphics[height=150mm,width=150mm]{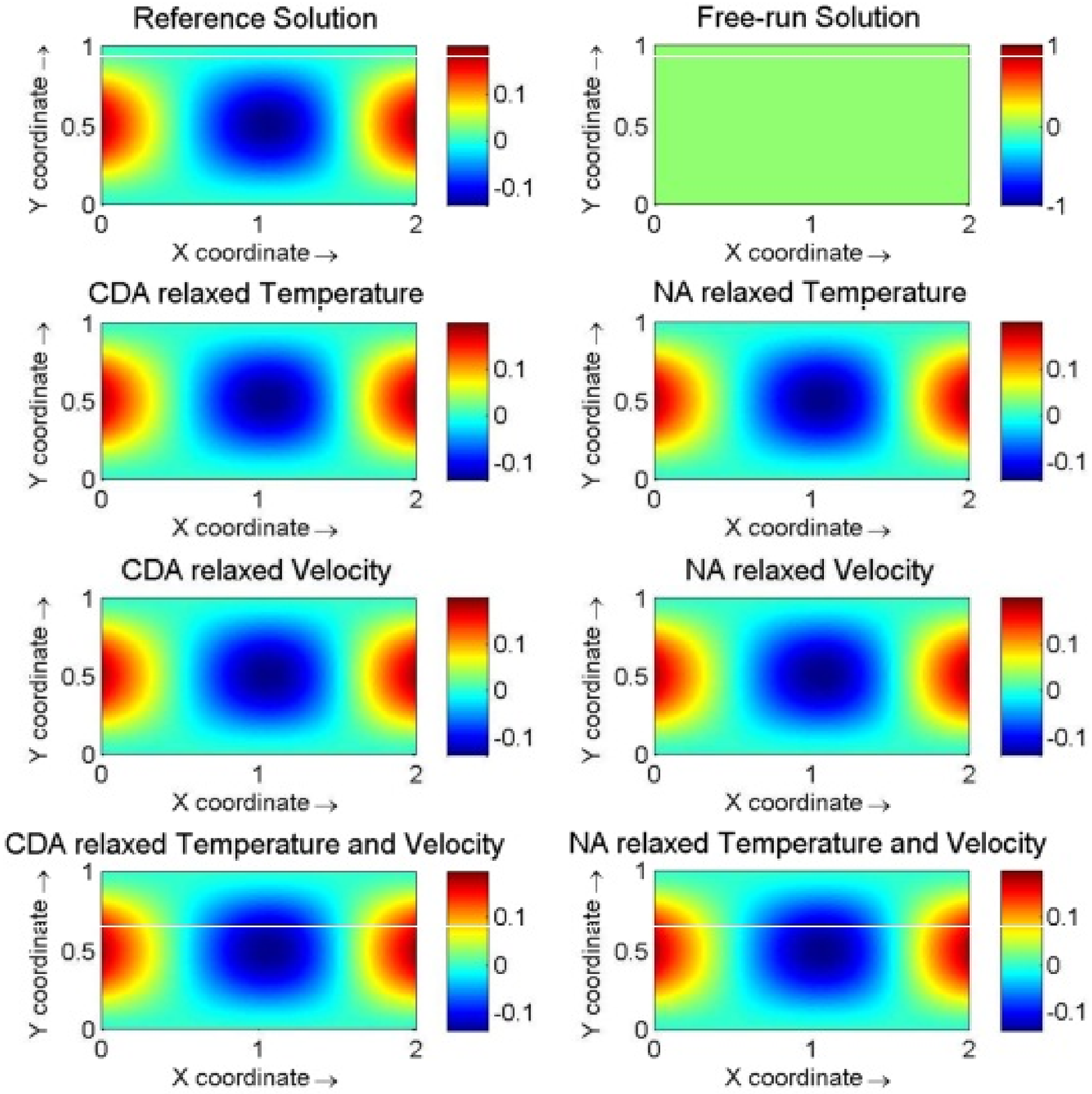}
\vskip 20pt 
	\caption{$y$-direction of velocity field at the end of the assimilation window ($t = 30$) and assimilation of data every 10 point grid resolution as it results from CDA and NA with assimilation of different variables (CDA - $\mu_\theta=0.10, \mu_V = 0.10$, NA - $\alpha_\theta=1.50, \alpha_V = 1.50$). The simulation starts from rest, i.e. $v=0$ (\ref{fig2} - left panel).}
\label{fig7}
\end{figure}

When assimilating only velocity (i.e., $\mu_\theta = 0$), both temperature and velocity converge, but velocity converges faster. In the case of velocity, the choice of relaxation parameter ($\mu_u$) plays an important role as can be seen from Figures \ref{fig2} - \ref{fig4}. Using larger values of $\mu_u$, the velocity converges to the reference solution exponentially. This suggests that for the given set of initial conditions, velocity and temperature measurements alone are enough to determine all the fields in the 2D B\'{e}nard convection model using CDA as demonstrated in \cite{Farhat2015}. The results further demonstrate that CDA is robust and not very sensitive to the choice of relaxation parameters ($\mu_\theta$ and $\mu_u$). The snapshots of the assimilated temperature and velocity data at the end of assimilation window are presented in Figures (\ref{fig5} - \ref{fig7}). These figures show that both CDA and NA converge to the reference solutions at the end of assimilation window.

\subsection{Sensitivity with respect to time frequency of the data}

We also investigated the sensitivity of CDA and NA to the frequency of available data in time. Here we assume that the data are available every $10$ and $20$ time steps and assimilated both temperature and velocity. We implemented four different strategies:
\begin{itemize}
\item NA-Time: NA with relaxation of the model solution only at the time when the observations are available every $10$ and $20$ time steps.
\item CDA-Time: CDA with relaxation of the model solution only at the time when the observations are available every $10$ and $20$ time steps.
\item NA: NA with relaxation of the model solution when the observations are available every time step.
\item CDA: CDA with relaxation of the model solution when the observations are available every time step.
\end{itemize}

\begin{figure}[htbp]
\centering
	\includegraphics[height=150mm,width=150mm]{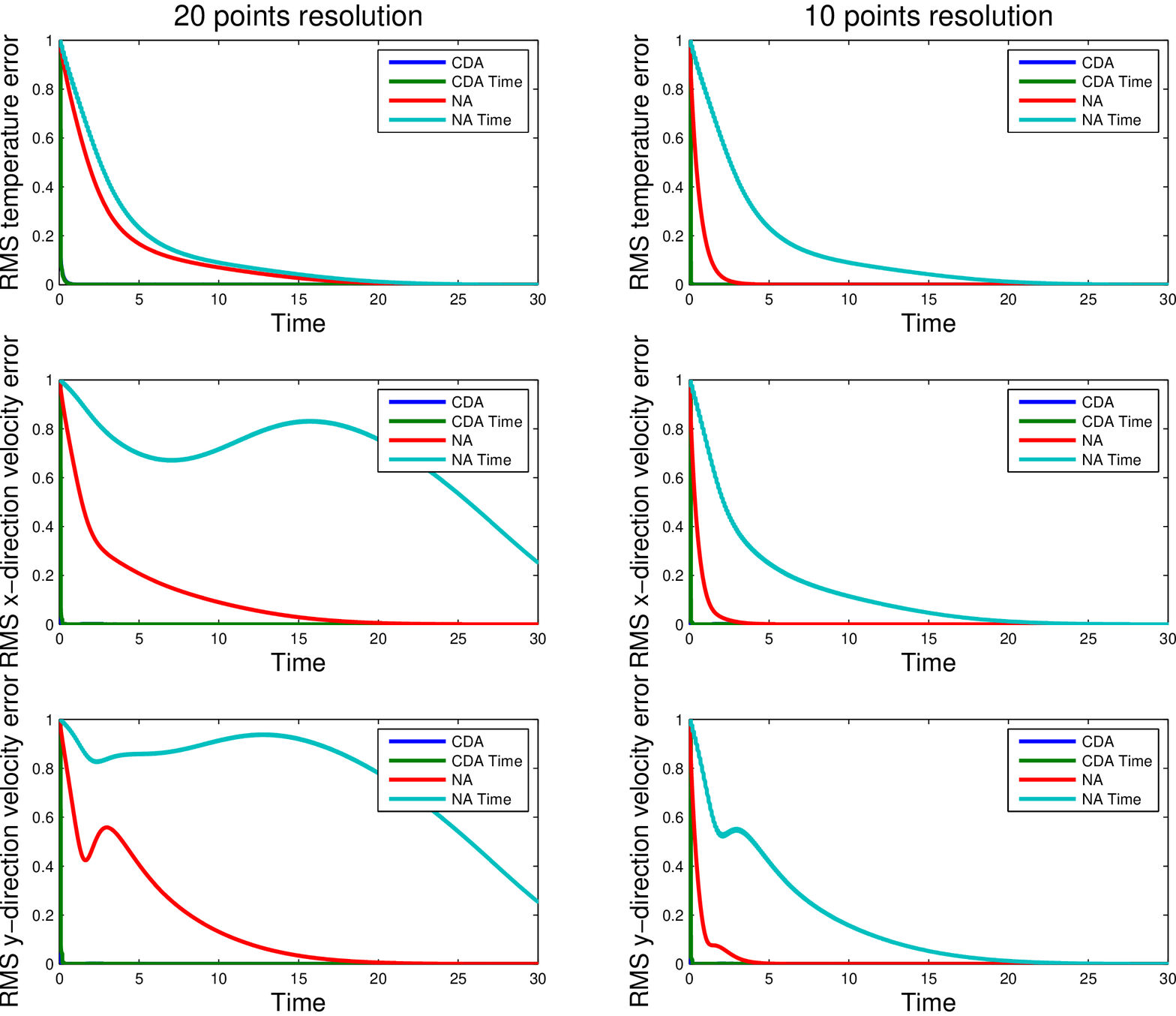}
\vskip 70pt 
	\caption{Sensitivity of CDA and NA to the availability of data every $10^{th}$ time steps (CDA - $\mu_\theta=0.50, \mu_V = 0.50$, NA - $\alpha_\theta=2.50, \alpha_V = 2.50$). The simulation starts from zero initial conditions (\ref{fig2} - left panel). The blue curve (CDA) is overlapped by the green line (CDA Time).}
\label{fig8}
\end{figure}

\begin{figure}[htbp]
\centering
	\includegraphics[height=150mm,width=150mm]{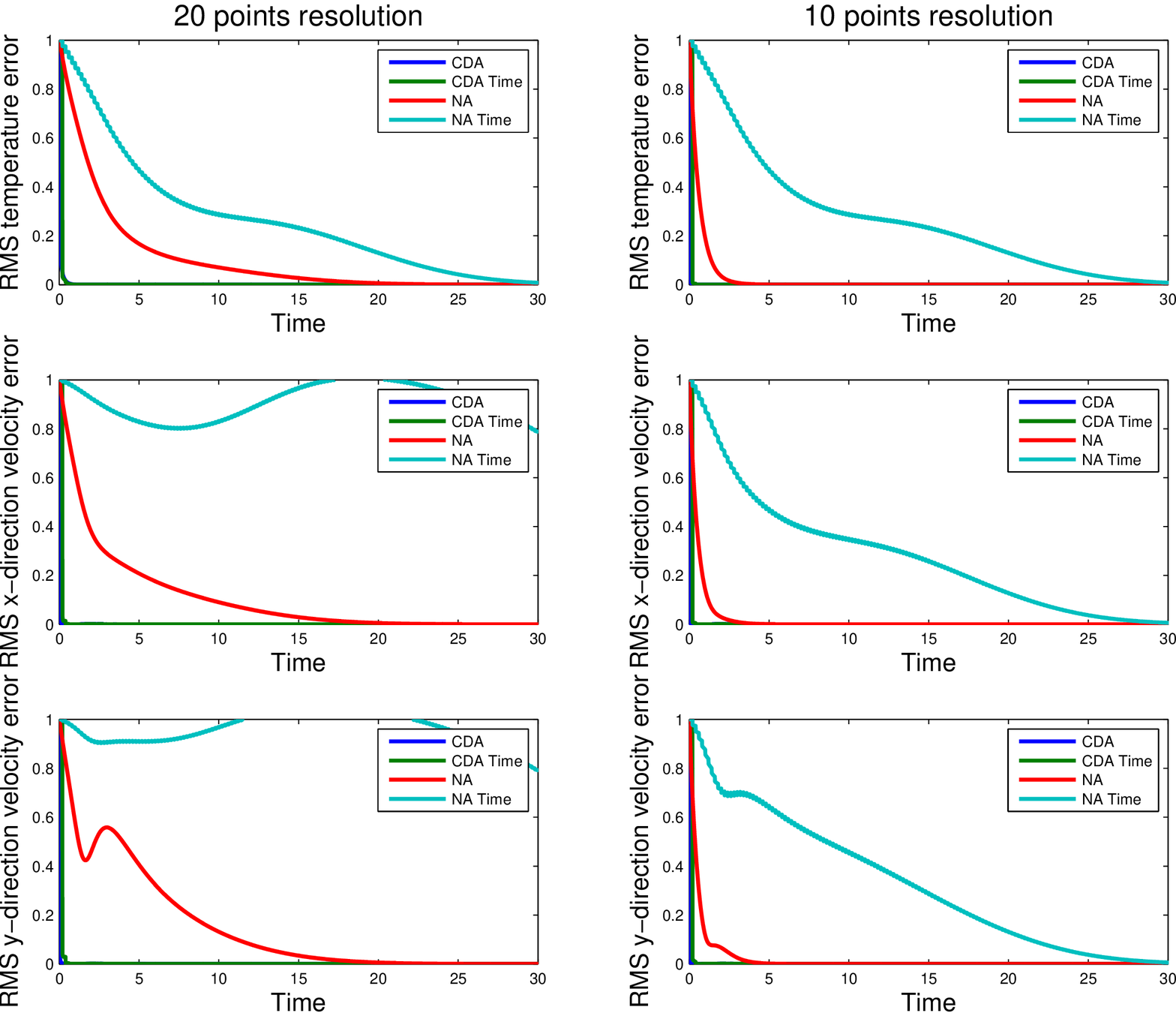}
\vskip 70pt 
	\caption{Sensitivity of CDA and NA to the availability of data every $20^{th}$ time steps (CDA - $\mu_\theta=0.50, \mu_V = 0.50$, NA - $\alpha_\theta=2.50, \alpha_V = 2.50$). The simulation starts from zero initial conditions (\ref{fig2} - left panel). The blue curve (CDA) is overlapped by the green line (CDA Time).}
\label{fig9}
\end{figure}

Figures \ref{fig8} - \ref{fig9} present the results of these experiments. The left and right panels of Figures \ref{fig8} - \ref{fig9} show the RRMSE with 20 and 10 point grid resolution respectively. Even when the data are assimilated every $10$ and $20$ time steps, CDA and CDA-time converge exponentially in time and there are no obvious differences in their RRMSEs. In contrast, the results of NA and NA-time significantly differ in terms of time and grid nudging.  This indicates that CDA is much more robust than NA in the case when the observations are coarse in time, and that it achieves fast convergence for both velocity and temperature.

\subsection{Sensitivity of CDA to the interpolation operator}

\begin{figure}[htbp]
\centering
	\includegraphics[height=150mm,width=150mm]{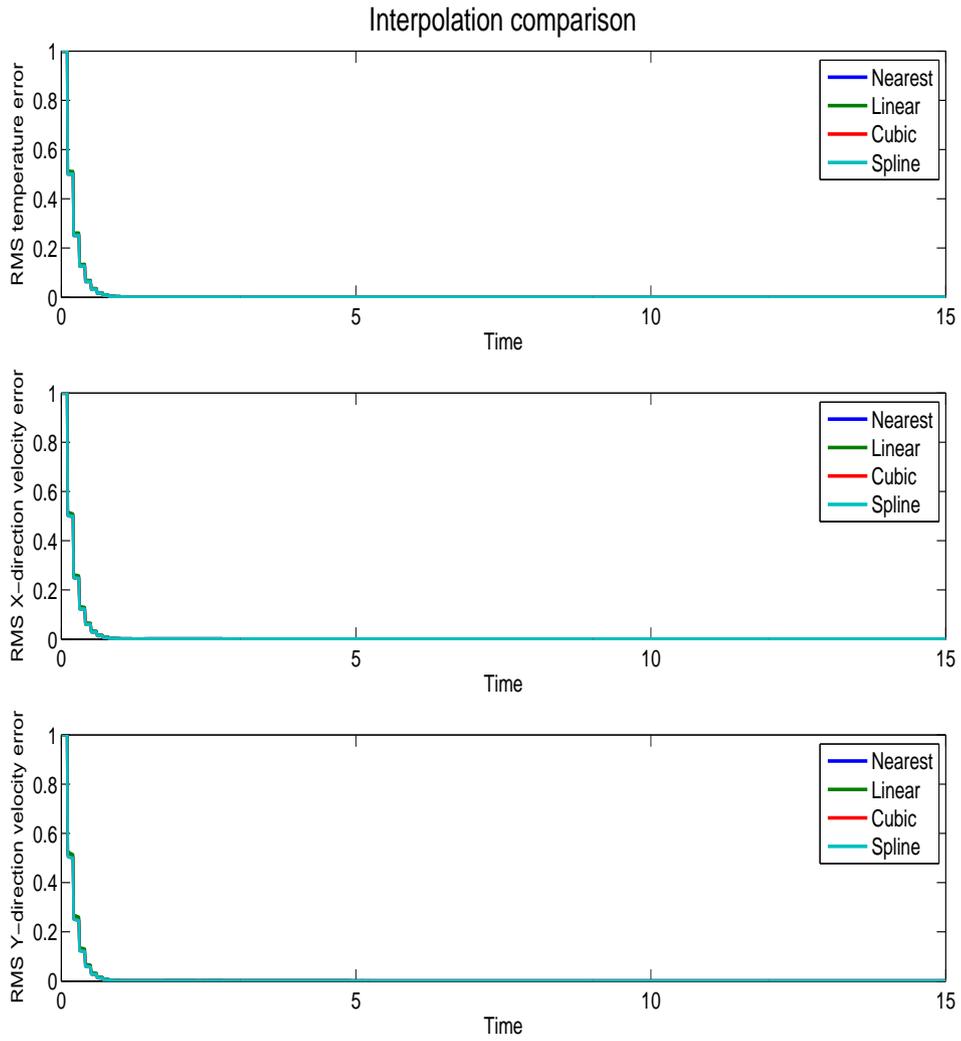}
\vskip 50pt 
	\caption{Sensitivity of CDA to the choice of the interpolation operator with $10$ point grid resolution every $10^{th}$ time step (CDA - $\mu_\theta=0.50, \mu_V = 0.50$, NA - $\alpha_\theta=2.50, \alpha_V = 2.50$).  The simulation starts from zero initial conditions (\ref{fig2} - left panel). All curves overlap and show similar convergence for all the interpolation operators.}
\label{fig10}
\end{figure}

We also study the sensitivity of CDA to the choice of the interpolation operator, examining  the impact of smoother interpolants. To do that, we set up an experiment in which the data are available $10$ time steps and $10$ grid points. Four different interpolation operators are tested (see Table \ref{tab3} for details).  The RRMSE results presented in Figure \ref{fig10} suggest that the CDA method is equally efficient for all the interpolation operators.

\begin{figure}[htbp]
\centering
	\includegraphics[height=150mm,width=150mm]{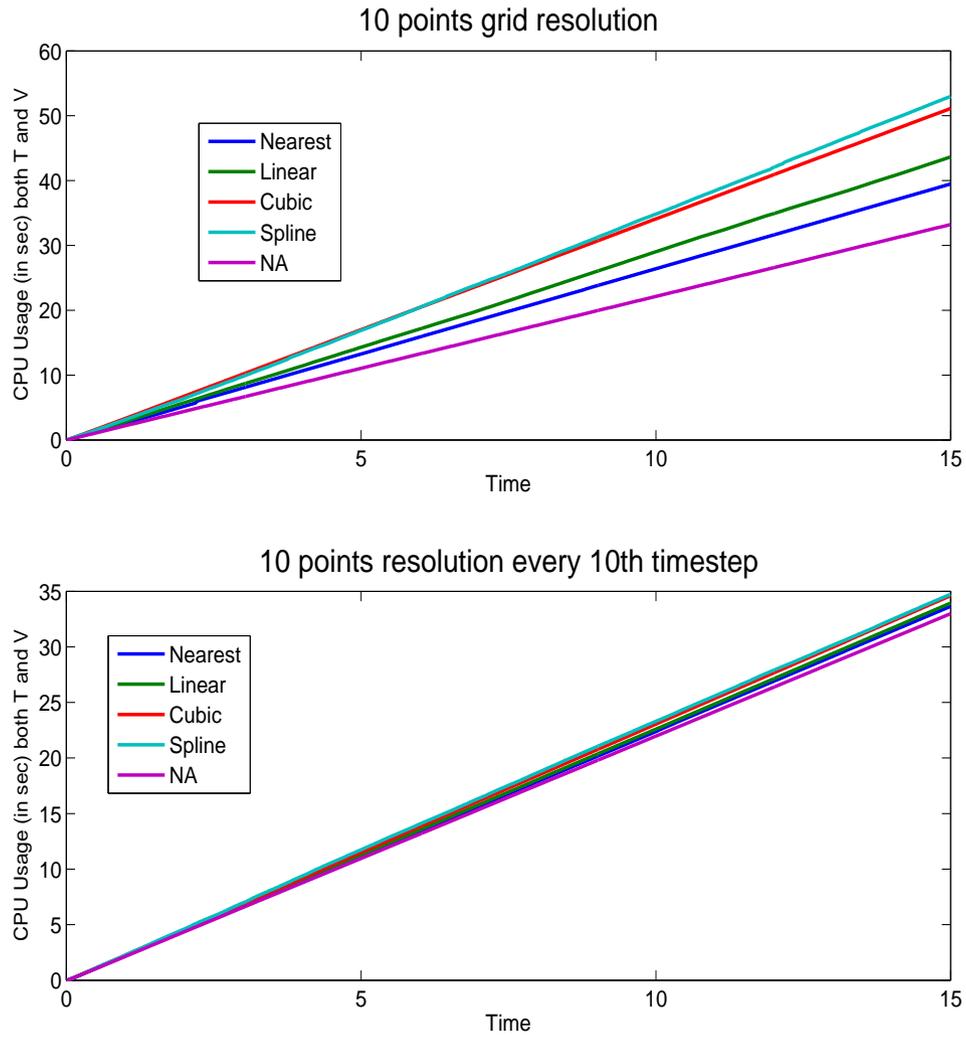}
\vskip 70pt 
	\caption{Efficiency of CDA vs NA for 10 point grid resolution (CDA - $\mu_\theta=0.50, \mu_V = 0.50$, NA - $\alpha_\theta=2.50, \alpha_V = 2.50$). The simulation starts from zero initial conditions (\ref{fig2} - left panel)}
	\label{fig11}
\end{figure}

In terms of computational cost, we examined the CDA method with respect to each interpolation operator. For this we have assumed that the data is available every $10^{th}$ grid point for 1) each time step and 2) every $10^{th}$ time step and assimilated both temperature and velocity. The costs of these runs are then compared to that of NA. Figure \ref{fig11} plots the CPU usage of CDA interpolation schemes vs NA. It is evident that CDA converges exponentially and there is no significant difference in terms of  computational costs of CDA with the different interpolation schemes compared to NA, especially when data are assimilated every $10$ time steps. 

\subsection{Sensitivity to observational errors}

The above experiments assimilated perfect observations, meaning that they were exactly the values as extracted from the reference model run. In practical applications, however, observations are collected by measurements and those are prone to noise. Therefore, it is important to investigate the sensitivity of CDA to data noise. 

\begin{figure}[htbp]
\centering
	\includegraphics[height=150mm,width=150mm]{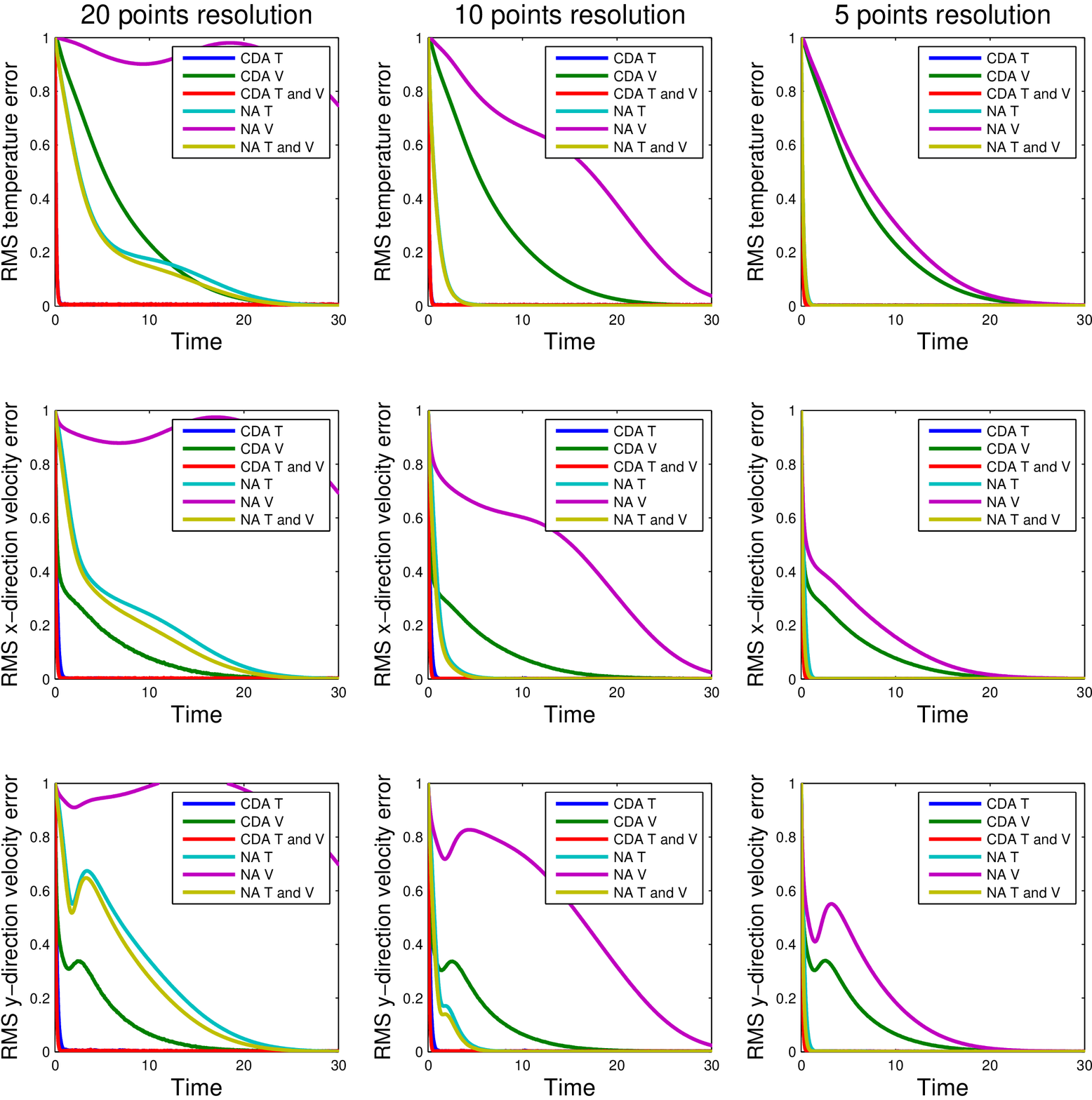}
\vskip 100pt 
	\caption{CDA  vs NA RMSE results with perturbed data ($5\%$) (CDA - $\mu_\theta=0.1, \mu_V = 0.1$, NA - $\alpha_\theta=1.50, \alpha_V = 1.50$). The left panel shows 20 point resolution, the middle panel shows 10 point resolution, and  5 point resolution on the right panel. The simulation starts from zero initial conditions (\ref{fig2} - left panel). The blue curve (CDA T) is overlapped by the red line (CDA T and V).}
\label{fig12}
\end{figure}

\begin{figure}[htbp]
\centering
	\includegraphics[height=150mm,width=150mm]{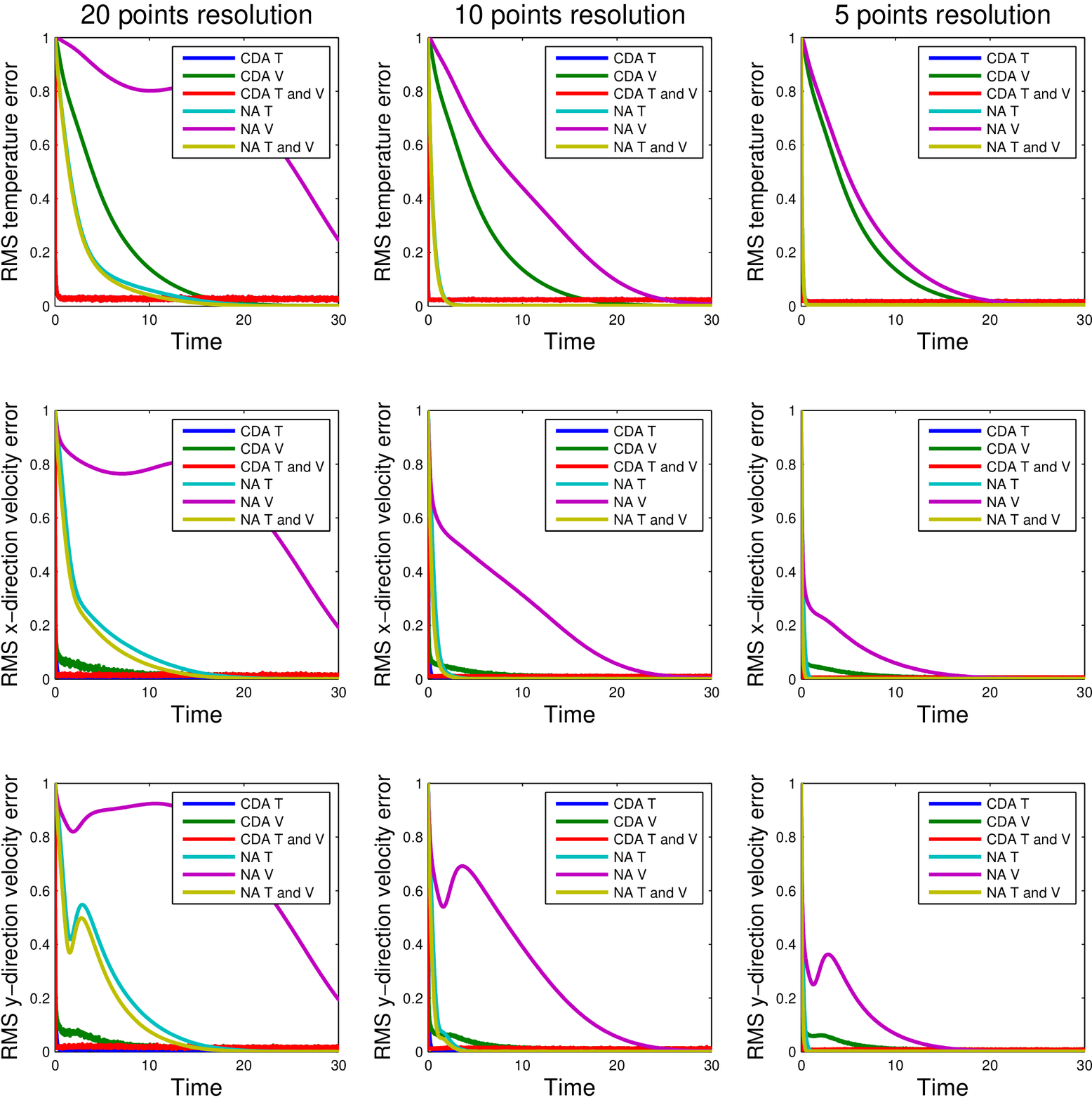}
\vskip 100pt 
	\caption{CDA  vs NA RMSE results with perturbed data ($5\%$) (CDA - $\mu_\theta=1.0, \mu_V = 1.0$, NA - $\alpha_\theta=3.50, \alpha_V = 3.50$). The left panel shows 20 point resolution, the middle panel shows 10 point resolution, and  5 point resolution on the right panel. The simulation starts from zero initial conditions (\ref{fig2} - left panel). The blue curve (CDA T) is overlapped by the red line (CDA T and V).}
\label{fig13}
\end{figure}

To simulate the situation of noisy observations, perturbations were induced in the process of obtaining the reference solution. Each entry in the observations set was multiplied by a random number between $1+\varepsilon$ and $1-\varepsilon$, corresponding to measuring with a relative error of $\varepsilon$. The experiments were performed with $\varepsilon$ values of 5\%. We have used two values of relaxation parameters ($\mu_\theta$ and $\mu_u$) and nudging operators ($\alpha_\theta$ and $\alpha_u$), referred to as large and small in Table \ref{tab2}. The data is assimilated every time step from $20$, $10$ and $5$ grid point resolution. The results are presented in Figures \ref{fig12} - \ref{fig13}. The RRMSE results are not impacted by the observational noise for small relaxation parameters ($\mu_\theta=0.10$ and $\mu_u=0.10$), in contrast to the CDA runs with large values of relaxation parameters ($\mu_\theta=1.0$ and $\mu_u=1.0$) which reveal a slight impact on the RRMSE values specially for the case of coarse data assimilating temperature and velocity every 20 point grid. This is expected since larger values of relaxation parameters nudge more strongly the model toward the perturbed observations. Also note that although NA exhibits slow convergence rate, its results are also sensitive to noise in the observations.

\subsection{Assimilation using temperature only: A counter example}

For the present setting, the possibility of assimilating all field variables based on observations of the velocity alone was established rigorously.
The numerical experiments above exhibited results that are consistent with theory.  The simulations also showed that in many cases,
coarse temperature observations can be sufficient for predicting both the temperature and velocity fields.  However, we do not expect
coarse temperature observations to be sufficient in general.  Intuitively, two scenarios immediately come to mind that signal potential deficiencies.
One concerns the situation where the Prandtl number is small.  In this case, temperature gradients tend to equilibrate at much
higher rate than regions of high shear, leading to insufficient information about the structure of the velocity field.  Another situation concerns the
existence of purely horizontal shear, which can exist in the absence of a temperature gradient.  

\begin{figure}[htbp]
\centering
	\includegraphics[height=150mm,width=150mm]{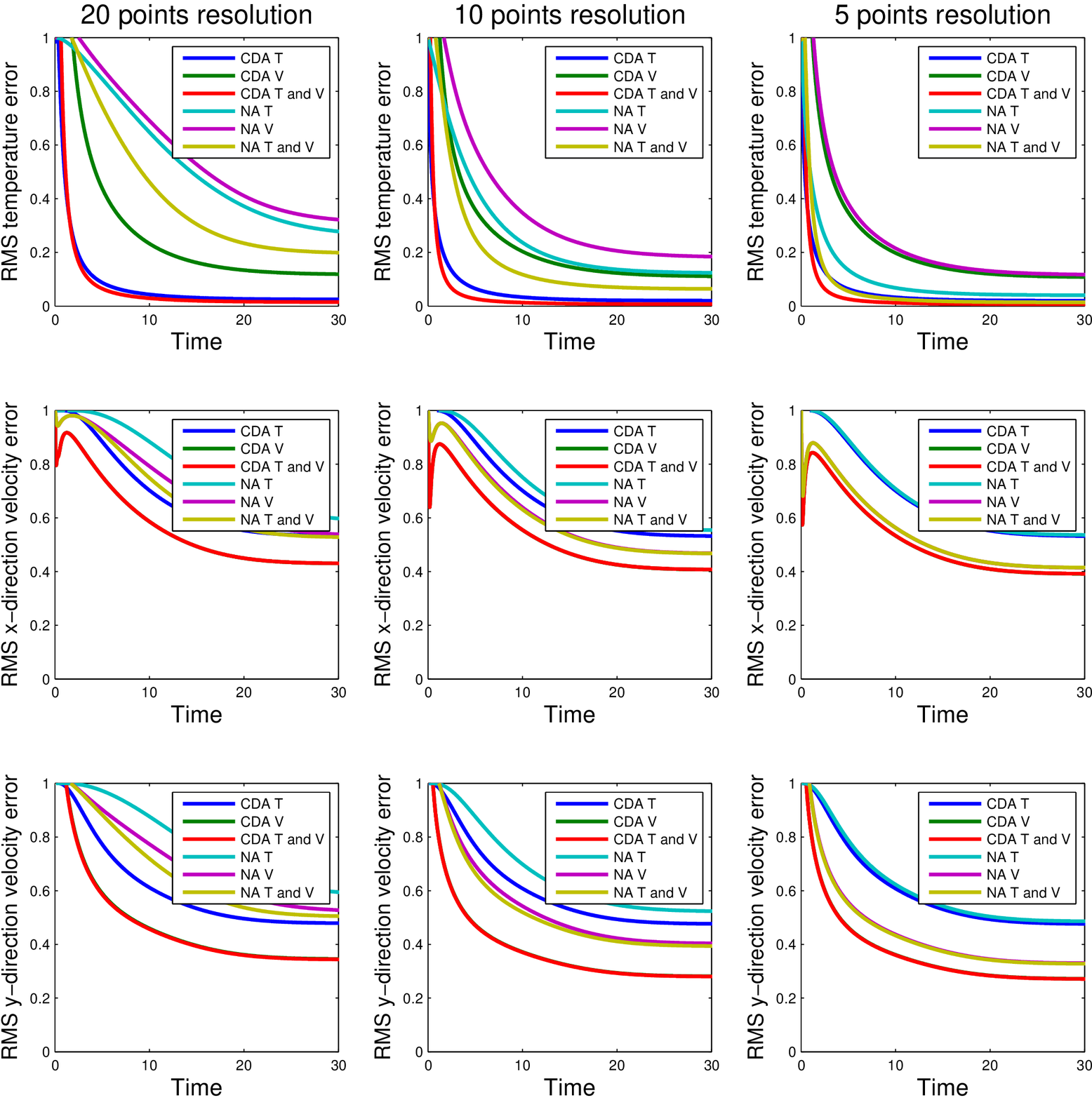}
\vskip 100pt 
	\caption{CDA  vs NA RMSE results (CDA - $\mu_\theta=0.1, \mu_V = 0.1$, NA - $\alpha_\theta=1.50, \alpha_V = 1.50$) with new set of initial condition. The left panel shows 20 point resolution, the middle panel shows 10 point resolution, and  5 point resolution on the right panel. }
\label{fig17}
\end{figure}

\begin{figure}[htbp]
\centering
	\includegraphics[height=150mm,width=150mm]{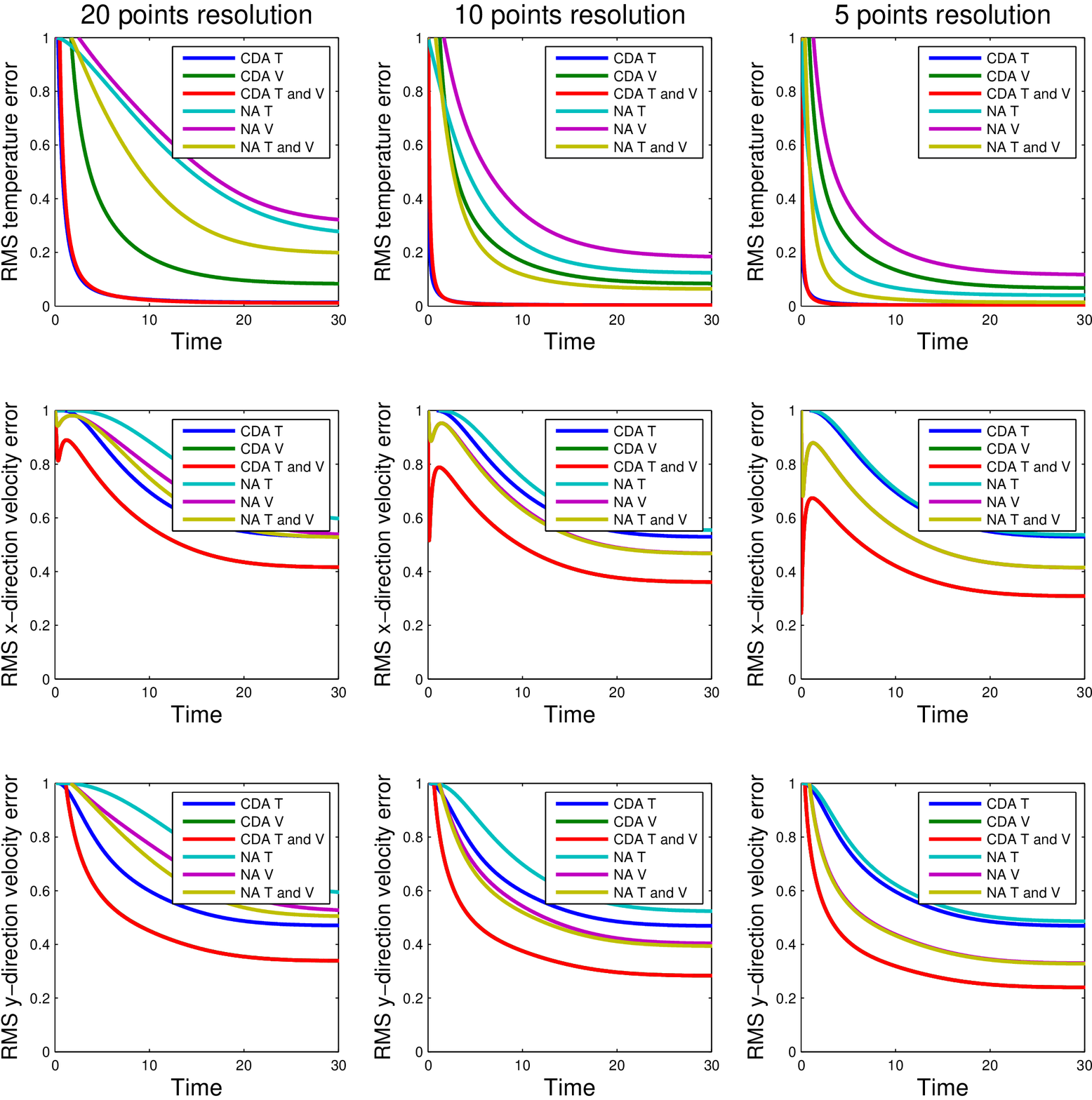}
\vskip 100pt 
	\caption{CDA  vs NA RMSE results (CDA - $\mu_\theta=1.0, \mu_V = 1.0$, NA - $\alpha_\theta=3.50, \alpha_V = 3.50$) with new set of initial condition. The left panel shows 20 point resolution, the middle panel shows 10 point resolution, and  5 point resolution on the right panel.}
\label{fig18}
\end{figure}

In this section, we conduct a numerical experiment, inspired by the latter observation, demonstrating that assimilation of temperature measurements 
only may not always recover the true state.  In the present experiment, the initial condition for the $x$-direction velocity is:
$$
u(x,y,0) = \left\{ \begin{array}{ll}
0.5sin(2 \pi / 20 y) & \mbox{for $0.4 \leq y \leq 0.6$}\\
0 & \mbox{otherwise} 
\end{array}
\right.
$$
representing a strong velocity shear in the center of the domain. Temperature and velocity in $y$ direction are set to zero (see Figure \ref{fig1}, right panel). The rest of the experimental setup is the same as in Section \ref{sengrid}. We performed three experiments assimilating $(i)$ velocity, $(ii)$ temperature, and $(iii)$ both velocity and temperature, with $20$, $10$ and $5$ point grid resolution. Figures \ref{fig17} - \ref{fig18}, compare the RRMSE with respect to the norm of the reference solution as they result from CDA and NA assimilating different types of data using two values of relaxation parameters (referred to as small and large in Table \ref{tab2}). The following conclusions can be drawn from this experiment:

\begin{itemize}
\item When assimilating only temperature with CDA (i.e., $\mu_u = 0$), the temperature converges exponentially, but velocity does not converge. 
\vskip 10pt
\item When assimilating velocity (i.e., $\mu_\theta = 0$) or both temperature and velocity with CDA, the temperature converges exponentially and velocity also converged. 
\vskip 10pt
\item Larger values of relaxation parameters, lead to lower RRMSE values in case of velocity. 
\end{itemize}

\begin{figure}[htbp]
\centering
	\includegraphics[height=150mm,width=150mm]{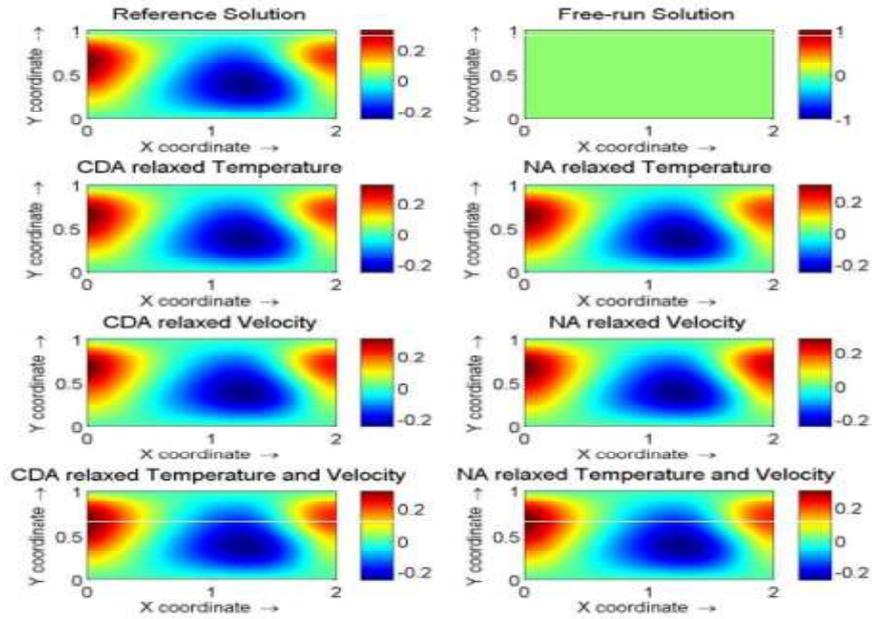}
\vskip 20pt 
	\caption{Temperature field at time ($t = 30$) and assimilation of data every 5 point grid resolution as it results from CDA and NA with new set of initial conditions (CDA - $\mu_\theta=1.0, \mu_V = 1.0$, NA - $\alpha_\theta=3.50, \alpha_V = 3.50$). The simulation starts from zero initial conditions (\ref{fig2} - left panel).}
\label{fig19}
\end{figure}

\begin{figure}[htbp]
\centering
	\includegraphics[height=150mm,width=150mm]{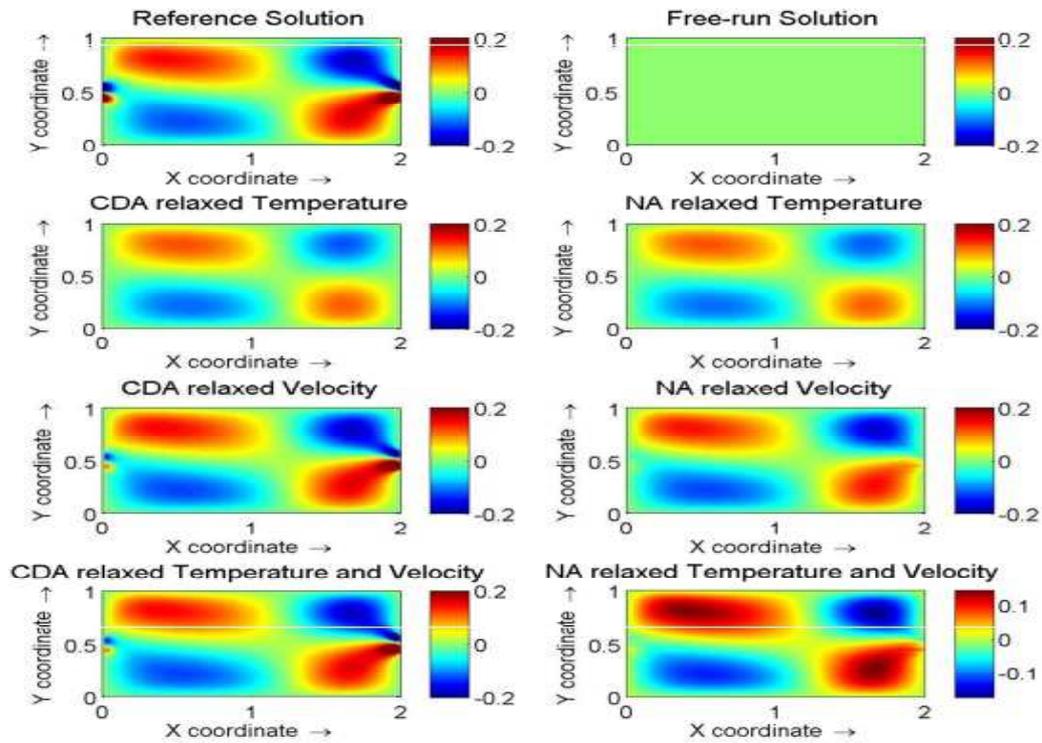}
\vskip 20pt 
	\caption{$x$-direction of velocity field at time ($t = 30$) and assimilation of data every 5 point grid resolution as it results from CDA and NA with with new set of initial conditions (CDA - $\mu_\theta=1.0, \mu_V = 1.0$, NA - $\alpha_\theta=3.50, \alpha_V = 3.50$). The simulation starts from rest, i.e., $u=0$ (\ref{fig2} - left panel).}
\label{fig20}
\end{figure}

\begin{figure}[htbp]
\centering
	\includegraphics[height=150mm,width=150mm]{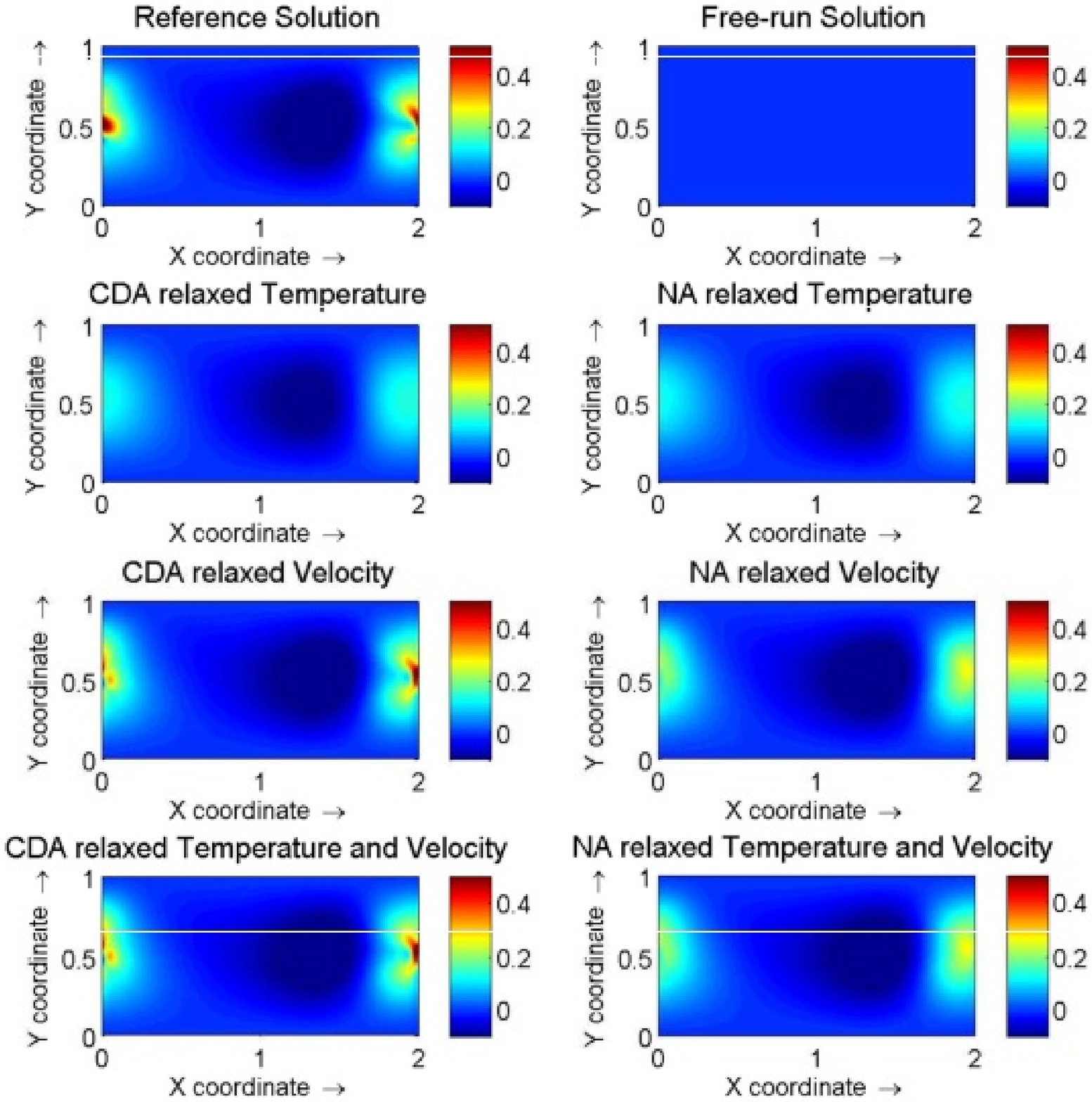}
\vskip 20pt 
	\caption{Y-direction of velocity field at the end of the assimilation window ($t = 30$) and assimilation of data every 5 point grid resolution as it results from CDA and NA with assimilation of different variables (CDA - $\mu_\theta=1.0, \mu_V = 1.0$, NA - $\alpha_\theta=3.50, \alpha_V = 3.50$). The simulation starts from rest, i.e., $v=0$ (\ref{fig2} - left panel). }
\label{fig21}
\end{figure}

The snapshots of the assimilated temperature and velocity data at the end of the assimilation window are presented in Figures \ref{fig19} - \ref{fig21}. These figures compare the assimilated results from CDA and NA to the reference solutions at the end of assimilation window. Consistent with the RRMSE results, one can clearly see that in this case the temperature measurements alone are not able to recover velocity fields at the end of the assimilation window.

\section{Conclusions}\label{sec.con}
Dynamical Downscaling (DD) methods are needed to compute local scale weather and climate from the outputs of global circulation models (GCMs). Nudging based methods are often used for DD to force regional circulation model (RCM) simulations towards large-scale driving data. In this paper the continuous data assimilation (CDA) algorithm is implemented and tested in the context of DD for the 2D B\'{e}nard convection model using coarse mesh measurements of the temperature and velocity. 

Various numerical experiments are performed with CDA assimilating velocity, temperature, and both velocity and temperature. The results obtained from these experiments are compared with the standard point-to-point grid nudging method and demonstrate that CDA is more robust than NA, and is not sensitive to the choice of control parameters. It is also shown that assimilating only velocity, both temperature and velocity converge, though velocity converges faster. On the other hand assimilating only temperature does not necessarily recover the velocity field. This is demonstrated using an example where assimilation starts from an initial shear flow. In particular the results indicate that in this case the estimated velocity field at the end of assimilation window suffers from large errors.

In addition we have examined the sensitivity of CDA and grid nudging methods with respect to perturbations in the assimilated observations (i.e. the values that are available either as measurements or from numerical experiments) and the spatial and temporal frequencies of observation points. In both cases, CDA  was found to be more performant and robust than the standard grid nudging. CDA was found slightly sensitive to the choice of the relaxation parameter when perturbations are applied to the observations. Even though CDA exhibits fast convergence, the computational costs of CDA and grid nudging are very similar. Our overall results suggest that CDA is easy to implement, and provides an efficient and robust approach for dynamical downscaling of the 2D B\'{e}nard convection model.

\subsection{Acknowledgments} The work of E.S.T. was supported in part by the ONR grant N00014-15-1-2333 and the NSF grants DMS-1109640 and DMS-1109645.


\clearpage

\clearpage

\begin{table}
\centering
\begin{tabular}{|c|c|c|}
\hline
Interpolant & Relaxation parameter (CDA) & Nudging parameter (NA) \\ \hline
Large  & $\mu_\theta = 1.0, \mu_u = 1.0$ & $\alpha_\theta = 3.5 , \alpha_u = 3.5$ \\ \hline
Medium & $\mu_\theta = 0.50, \mu_u = 0.50$  & $\alpha_\theta = 2.5 , \alpha_u = 2.5$ \\ \hline
Small & $\mu_\theta = 0.10, \mu_u = 0.10$ & $\alpha_\theta = 1.5 , \alpha_u = 1.5$ \\  
\hline
\end{tabular}
\caption{Relaxation and nudging parameters for CDA and NA respectively.}
\label{tab2}
\end{table}

\clearpage

\begin{table}
\centering
\begin{tabular}{|c|c|c|}
\hline
Interpolant &  Details \\ \hline
Nearest  & Discontinuous \\ \hline
Linear & $C([0, 1] \times [0,1])$ \\ \hline
Cubic & $C^1([0, 1] \times [0,1])$ \\ \hline
Spline & $C^2([0, 1] \times [0,1])$ \\
\hline
\end{tabular}
\caption{2D Interpolation operators.}
\label{tab3}
\end{table}
\clearpage

\clearpage

\clearpage
\clearpage

\clearpage

\clearpage

\clearpage

\clearpage

\clearpage
\clearpage
\clearpage
\clearpage

\clearpage
\clearpage

\clearpage

\clearpage

\clearpage

\clearpage

\end{document}